\documentclass[preprint,authoryear,12pt]{elsarticle}

\usepackage{amssymb}

\newdefinition{rmk}{Remark}

\journal{ }

\begin{document}

\begin{frontmatter}



\title{Composition and inverse of multivariate functions and algebraic system of equations }


\author{Wu Zi qian}

\address{Email:humble175@sohu.com.Address:Shenzhen city,China.}

\begin{abstract}

 Every one knows that an equation is equivalent to a multivariate function. Generally speaking, there are more than one unknown x in this multivariate function and it is not easy to reduce the number of unknown x to one. In this paper we achieve this by introducing function promotion which can converse a function of less variables to one of more variables and by introducing multivariate function composition developed from unary function composition. We introduced inverse multivariate functions extended from inverse unary functions then we can express the solution by an inverse multivariate function if this equivalent multivariate function is invertible.  For an equivalent irreversible multivariate function we introduced relation and consider the equivalent multivariate function as a special multivariate relation then we can express the solution by an inverse multivariate relation which always exists. No one will belief that every one is familiar with multivariate function composition after their reading this paper and knowing the importance of  multivariate function composition to express the solution of a general equation. Further more we discuss the possibility of expressing the solution, a multivariate function or multivariate relation, for a general equation by superposition of unary ones. This will involve Hilbert's 13th problem. The topic shown in this paper can be new directions for many mathematicians because the topic is so basic.

\end{abstract}


\end{frontmatter}
\section{Introduction}
\label{intro}

It is an embarrassing thing to talk about'The explicit solution for a general equation'. Nearly every mathematician hate it.  First we show our key view points and the strategy about 'How to give an explicit solution for a general equation?' then we mention briefly the history about this topic to avoid the argue.

For clarity we give new symbolics for several known operations, $f_{a}$ for addition,$x_{1}+x_{2}=f_{a}(x_{1},x_{2})$, $f_{s}$ for subtraction, ,$x_{1}- x_{2}=f_{s}(x_{1},x_{2})$, $f_{m}$ for multiplication, $x_{1}\times x_{2}=f_{m}(x_{1},x_{2})$, $f_{d}$ for division, $x_{1}\div x_{2}=f_{d}(x_{1},x_{2})$,  $f_{p}$ for power,  $x_{1}^{x_{2}}=f_{p}(x_{1},x_{2})$, $f_{r}$ for root,  $\sqrt[x_{2}]{x_{1}}=f_{r}(x_{1},x_{2})$, $f_{l}$ for logarithm,$\log_{x_{1}} x_{2}=f_{l}(x_{1},x_{2})$, respectively.

   First, an equation is equivalent to a multivariate function. For example, for the equation
$x^{a}+x^{b}=c$ or $f_{a}[f_{p}(x,a),f_{p}(x,b)]=c$, the left of it is equivalent to a multivariate function f(x,a,b). In this paper, we always left a known parameter on the right of =, like c in this equation.  However, f(x,a,b) is a meaningless expression for it dealing with  no any element of the equation,$f_{a}$ or $f_{p}$. Our strategy is giving f a expression like $f=w(f_{p},f_{a})$,that is say f is a function of $f_{p}$ and $f_{a}$ or f is constructed by $f_{p}$ and $f_{a}$. We achieve this point by multivariate function composition which are developed from unary function composition and  function promotion which can converse a function of less variables to one of more variables.

We extended inverse unary functions to inverse multivariate functions then we can express the solution by an inverse multivariate function if this equivalent multivariate function is invertible.  For an equivalent irreversible multivariate function we introduced relation and consider the equivalent multivariate function as a special multivariate relation then can express the solution by an inverse multivariate relation which always exists.

Solution expressed by a multivariate function or a multivariate relation is called  multivariate solution. No one doubt its existing. This paper shown a perfect way to gave multivariate solution.

Can we transform these multivariate solution into binary ones or unary ones? This involves whether or not we can express the  multivariate function by binary functions or superpositions of unary functions. Hilbert tended to a negative answer and wrote his doubt into his famous 23 problems, the 13th problem[3]. A.N.Kolmogorov gave a conclusion that any multivariate continuous functions can be expressed by unary continuous functions[2].There is argue about expressing multivariate continuous functions  by superposition of unary analytic functions or by unary smooth functions[1].

By the result above we can express  multivariate solutions by unary continuous functions or unary relations. However these continuous functions are at best only continuous.

In special situation we can gave solutions expressed by binary smooth functions or binary analytic functions. This will be an important direction for us.

There are special equations constructed by functions over finite sets. We call these functions discrete functions and these equations discrete equations. Then we shown a conclusion that any multivariate discrete function can always be expressed by superposition of unary discrete functions so discrete equations have binary solutions. Discrete equations are so easy and many new concepts and new ways will be developed from them.

\section{Core concepts}

 We will show several new concepts in this section.

{\bf ~2.1 Multivariate function composition and inverse multivariate functions}

{\bf Definition~2.1 Multivariate function composition}  Let B be a set, $B\neq \emptyset$, $F^{n}:=\{f|f:B^{n}\longrightarrow B\}$,  $\mathbf{C_{i}}$ is defined as:

$\mathbf{C}_{i}:F^{n}\times F^{n}\longrightarrow F^{n},(f_{1},f_{2})\longmapsto f_{3}$ and

 $f_{3}(x_{1},\ldots,x_{n})=f_{1}[x_{1},x_{2},\cdots,x_{i-1},f_{2}(x_{1},\ldots,x_{n}),x_{i+1},\cdots,x_{n}]$

We denote it as:

\begin{equation}\label{eq:eps}
f_{3}=\mathbf{C}_{i}(f_{1},f_{2})
\end{equation}

or

\begin{equation}\label{eq:eps}
f_{3}=\mathbf{C}_{i}\frac{f_{1}}{f_{2}}
\end{equation}

We ought to think of $\mathbf{C}_{i}-$ as an entirety and do not misunderstander $\frac{f_{1}}{f_{2}}$ as a fraction.

Aim of introducing new concepts is to solve equations so the set B  will be limited in N,Z,Q,R,C or their subsets in this paper .

$f_{1}$ or $f_{2}$ may lack several variables.

\textbf{Example 1} $f_{1}$ is a function of n variables,$f_{1}=g[x_{1},x_{2},\cdots,x_{i-1},x_{i},x_{i+1},\ldots,x_{n}]$, $f_{2}$ is an unary function, $f_{2}=\beta$,

$f_{3}[x_{1},x_{2},\cdots,x_{i-1},x_{i},x_{i+1},\ldots,x_{n}]=[\mathbf{C}_{i}\frac{f_{1}}{\beta}][x_{1},x_{2},\cdots,x_{i-1},x_{i},x_{i+1},\ldots,x_{n}]$

$=g[x_{1},x_{2},\cdots,x_{i-1},\beta(x_{i}),x_{i+1},\ldots,x_{n}]$��

\textbf{Example 2}  $f_{1}$ is a function of n variables,$f_{1}=g[x_{1},x_{2},\cdots,x_{i-1},x_{i},x_{i+1},\ldots,x_{n}]$, $f_{2}$ is a constant,$f_{2}=a$��

$f_{3}[x_{1},x_{2},\cdots,x_{i-1},x_{i},x_{i+1},\ldots,x_{n}]=[\mathbf{C}_{i}\frac{f_{1}}{a}][x_{1},x_{2},\cdots,x_{i-1},x_{i+1},\ldots,x_{n}]$

$=f_{1}[x_{1},x_{2},\cdots,x_{i-1},c,x_{i+1},\ldots,x_{n}]$

$f_{3}$ is a function of n-1 variables, a special function of n variables.

\textbf{Example 3} $f_{1}$ is an unary function, $f_{1}=\beta$, $f_{2}$ is a function of n variables,$f_{2}=g[x_{1},x_{2},\cdots,x_{i-1},x_{i},x_{i+1},\ldots,x_{n}]$,

$f_{3}[x_{1},x_{2},\cdots,x_{i-1},x_{i},x_{i+1},\ldots,x_{n}]=[\mathbf{C}_{i}\frac{\beta}{f_{1}}][x_{1},x_{2},\cdots,x_{i-1},x_{i},x_{i+1},\ldots,x_{n}]$

$=\beta[g(x_{1},x_{2},\cdots,x_{i-1},x_{i},x_{i+1},\ldots,x_{n})]$��

\textbf{Example 4} $f_{1}$ is a constant, $f_{1}=c$, $f_{2}$ is a function of n variables,$f_{2}=g[x_{1},x_{2},\cdots,x_{i-1},x_{i},x_{i+1},\ldots,x_{n}]$,there will not be substitution and $f_{3}=c$��

\textbf{Example 5} $f_{1}$ and $f_{2}$ are unary functions, $f_{1}=\beta_{1}$, $f_{2}=\beta_{2}$, $f_{3}=\mathbf{C}(f_{1},f_{2})=\beta_{1}.\beta_{2}$

{\bf Definition~2.2 Function promotion}  Let B be a set, $B\neq \emptyset$, $F^{n}:=\{f|f:B^{n}\longrightarrow B\}$,  $F^{m}:=\{f|f:B^{m}\longrightarrow B\}$,

Function promotion $\mathbf{P^{n}_{j_{1},j_{2},\cdots,j_{m}}}$ is defined as: $\mathbf{P}^{n}_{j_{1},j_{2},\cdots,j_{m}}:F^{m}\longrightarrow F^{n},g\longmapsto f$.

$f(x_{1},\ldots,x_{n})=g(x_{j_{1}},x_{j_{2}},\cdots,x_{j_{m}})+\cdots+O(x_{j_{m+1}})+\cdots+O(x_{j_{n}})$ where $1\leq j_{k}\leq n$ and $O(x_{j_{k}})\equiv0 $ for $k\geq m+1$

Only $x_{j_{1}},x_{j_{2}},\cdots,x_{j_{m}}$ are really existing variables, others should be omitted.

We denote it as:

\begin{equation}\label{eq:eps}
f=\mathbf{P}^{n}_{j_{1},j_{2},\cdots,j_{m}}(g)
\end{equation}

Usually we construct the expressions by binary operations or binary functions and we consider these binary functions as special multivariate functions then we ought to converse binary functions to multivariate functions by promotion frequently.

\textbf{Example 6} $x^{a}+x^{b}$, let $x_{1}=x,x_{2}=a,x_{3}=b$

$w_{1}(x,a,b)=x^{a}\Longleftrightarrow w_{1}=\mathbf{P}^{3}_{1,2}(f_{p})$, $x^{a}=\mathbf{P}^{3}_{1,2}(f_{p})(x,a,b)$

$w_{2}(x,a,b)=x^{b}\Longleftrightarrow w_{2}=\mathbf{P}^{3}_{1,3}(f_{p})$, $x^{b}=\mathbf{P}^{3}_{1,3}(f_{p})(x,a,b)$

$w_{3}(x,a,b)=a+b\Longleftrightarrow w_{3}=\mathbf{P}^{3}_{2,3}(f_{a})$, $a+b=\mathbf{P}^{3}_{2,3}(f_{a})(x,a,b)$

The same binary function $f_{p}$ is conversered to different functions of 3 variables $w_{1}$ and $w_{2}$ respectively.

By substituting $\mathbf{P}^{3}_{1,2}(f_{p})(x,a,b)$ for a of $a+b=\mathbf{P}^{3}_{2,3}(f_{a})(x,a,b)$ we obtain: $x^{a}+b=\mathbf{C}_{2}(\mathbf{P}^{3}_{2,3}(f_{a}),\mathbf{P}^{3}_{1,2}(f_{p}))(x,a,b)$

By substituting $\mathbf{P}^{3}_{1,3}(f_{p})(x,a,b)$ for b of $\mathbf{C}_{2}(\mathbf{P}^{3}_{2,3}(f_{a}),\mathbf{P}^{3}_{1,2}(f_{p}))(x,a,b)+b$ we obtain: $x^{a}+x^{b}=[\mathbf{C}_{3}\frac{\mathbf{C}_{2}(\mathbf{P}^{3}_{2,3}(f_{a}),\mathbf{P}^{3}_{1,2}(f_{p}))}{\mathbf{P}^{3}_{1,3}(f_{p})}](x,a,b)$

Here we can see the necessity of writing  addition as $f_{a}$ or power as $f_{p}$. By function promotion we make
$x^{a}\Rightarrow \mathbf{P}^{3}_{1,2}(f_{p})(x,a,b)$, $x^{b}\Rightarrow \mathbf{P}^{3}_{1,3}(f_{p})(x,a,b)$ and $a+b \Rightarrow \mathbf{P}^{3}_{1,2}(f_{p})(x,a,b)$
then we obtain a final function of 3 variables $[\mathbf{C}_{3}\frac{\mathbf{C}_{2}(\mathbf{P}^{3}_{2,3}(f_{a}),\mathbf{P}^{3}_{1,2}(f_{p}))}{\mathbf{P}^{3}_{1,3}(f_{p})}](x,a,b)$ by
composition.

 Can we understander the huge difference between $[\mathbf{C}_{3}\frac{\mathbf{C}_{2}(\mathbf{P}^{3}_{2,3}(f_{a}),\mathbf{P}^{3}_{1,2}(f_{p}))}{\mathbf{P}^{3}_{1,3}(f_{p})}](x,a,b)$ and $x^{a}+x^{b}$. Look at! So deep the secret is hidden. We did not really know how binary functions construct the expressions of multivariate functions.

{\bf Definition~2.3 Inverse multivariate functions}  Let B be a set, $B\neq \emptyset$, $F^{n}:=\{f|f:B^{n}\longrightarrow B\}$,

Inverse multivariate function $f_{i}$ about $x_{i}$ of $f$ is defined as:  $x_{i}=f_{i}(x_{1},\ldots,$
$x_{i-1},x_{0},x_{i+1},\ldots,x_{n})$ if $x_{0}=f(x_{1},\ldots,x_{i-1},x_{i},x_{i+1},\ldots,x_{n})$ and for any $x_{j}(j=1,2,\cdots,n,j\neq i)$ $x_{i}\mapsto x_{0}$ is invertible.

Introduce  $\mathbf{I}_{i}:F^{n}\longrightarrow F^{n},f\mapsto f_{i}$, the relation is:

\begin{equation}\label{eq:eps}
f_{i}=\mathbf{I}_{i}(f)
\end{equation}

\textbf{Example 7} For binary operations,  $f_{s}=\mathbf{I}_{1}(f_{a})=\mathbf{I}_{2}(f_{a})$, $f_{d}=\mathbf{I}_{1}(f_{m})=\mathbf{I}_{2}(f_{m})$, $f_{r}=\mathbf{I}_{1}(f_{p})$, $f_{l}=\mathbf{I}_{2}(f_{p})$.

{\bf ~2.2 The multivariate relation and inverse multivariate relation}

We can not obtain the inverse of function 'f' if $x_{0}=f(x_{1},\ldots,x_{i-1},x_{i},x_{i+1},$ $\ldots,x_{n})$ is irreversible in the definition 2.3. Of course applied mathematicians introduced many-valued functions,however these informal functions can not be accepted by mathematicians in pure mathematics majors. For this reason we introduce a new concept, multivariate relation.

A multivariate relation is defined upon the concept of n+1-ary relation which has been introduced already in discrete  mathematics and was not limited in discrete set.

{\bf Definition~2.4 n+1-ary relation} A n+1-ary relation R over non-empty sets $ B_{i}(1\leq i\leq n)$ and $B_{0}$ is a subset of their Cartesian product written as $R\subset B_{1}\times B_{2}\ldots \times B_{n}\times B_{0}=\{\langle b_{1},b_{2},\ldots,b_{n},b_{0}\rangle|b_{i}\in B_{i}(1\leq i\leq n),b_{0}\in B_{0}\}$.
In other words, R is a n+1-ary relation if any of its elements is an ordered n+1 tuple $<x_{1},x_{2},\cdots x_{n},x_{0}>$,($x_{i}\in B_{i},i=0\sim n$).
Specially R is called n+1-ary relation over B if $R\subset \underbrace{B \times B\ldots \times B\times B}_{n+1}=B^{n+1}$.

{\bf Definition~2.5 Multivariate relation}  Let B is a  set, R is a n+1-ary relation over B, we shall call it an relation of n variables if we take  i-th element of its ordered n+1-tuple as its  i-th variable and the last element as its value. We denote it by $R$.

We use the form of function to denote relations of M variables:

 \qquad\qquad \qquad\qquad$x_{0}=R(x_{1},\cdots,x_{i-1},x_{i},x_{i+1},\cdots,x_{n})$

{\bf Definition~2.6 Inverse of a multivariate relation:} For a relation of n variables, its inverse about i-th variable $R_{i}$ is defined as:

\begin{scriptsize}
$R_{i}:=\{\langle x_{1},x_{2},\ldots ,x_{i-1},x_{0},x_{i+1},\ldots, x_{n},x_{i}\rangle | \langle x_{1},x_{2},\ldots ,x_{i-1},x_{i},x_{i+1},\ldots, x_{n},x_{0} \rangle \in R\}$
\end{scriptsize}

\begin{equation}\label{eq:eps}
x_{i}=R_{i}(x_{1},\cdots,x_{i-1},x_{0},x_{i+1},\cdots,x_{n}),
\end{equation}

We introduce inverse operator $\mathbf{I}_{i}$ and express $R_{i}$ by unary operator $\mathbf{I}_{i}$:

\begin{equation}\label{eq:eps}
 R_{i}=\mathbf{I}_{i}(R)
\end{equation}

{\bf Definition~2.7 Transform an irreversible multivariate function into a multivariate relation:} Any relation is always invertible. We can consider
an irreversible  multivariate function g as a special multivariate relation R. We say we transform an irreversible  multivariate function g into a multivariate relation R and introduce transform operator $\mathbf{T}$ to express R as:

\begin{equation}\label{eq:eps}
 R=\mathbf{T}(g)
\end{equation}

Multi-valued function has been used by mathematicians. We must admit that it may make confuse to pure mathematics, such as algebra and topology. However we defined multivariate relation to avoid this.  Actually we have to introduce multivariate relation. So many irreversible functions appear when we study equations   and we have to discuss their inverse. Surely mathematicians in algebra and topology need not to face equations!

{\bf ~2.3 Algebraic System of Equations}
   An algebraic system of equations (B,E) is constructed by a non-empty set B and a set E of equations constructed by the binary functions over B.
If elements of B are numbers then equations in (B,E) are algebraic equations. B can be a set of complex numbers or of real number or of integer. Algebraic systems of equations with a set of real number or with a set of complex numbers will be our main target. If elements of B are functions then equations in (B,E) are operator equations.

\section{Multivariate  solution}

{\bf 3.1 Examples of multivariate solution}

   We can obtain multivariate solution of any equations,algebraical equations or operator equations by above new concepts. We shown the main skills by following
examples. We shown solving several types of equations. The procedure includes 3 steps. Step 1, All binary operations or binary functions on the left of the
equation should be converted to functions of n variables if the equation contains n parameters including unknown x in its left. However not every variable really exists and we will point out which ones really exist. Step 2, The left of the equation should be conversed to a function of n variables  by composition then we can reduce the number of unknown x in the expression to one. Step 3, Place the only unknown x to the left of '=' and all others to the right by inverse function.

\textbf{Example 8} $f_{2}[x,f_{1}(x,a)]=b$

There is no step1 so begin from step 2,$[\mathbf{C}_{2}(f_{2},f_{1})](x,a)=b$

Step 3,$x=\mathbf{I}_{1}\Big(\mathbf{T}[\mathbf{C}_{2}(f_{2},f_{1})]\}\Big)(b,a)$

$f_{3}\{x,f_{2}[x,f_{1}(x,a)]\}=b$

Begin from step 2,$[\mathbf{C}_{2}\frac{f_{3}}{\mathbf{C}_{2}(f_{2},f_{1})}](x,a)=b$

Step 3,$x=\mathbf{I}_{1}\Big(\{\mathbf{T}[\mathbf{C}_{2}\frac{f_{3}}{\mathbf{C}_{2}(f_{2},f_{1})}]\}\Big)(b,a)$

$f_{n}\cdots f_{3}\{x,f_{2}[x,f_{1}(x,a)]\}=b$ has a binary solution.

In the equation  $f_{2}[x,f_{1}(x,a)]=b$, let $f_{1}(x,a)=x$, $f_{2}=f$,b=a, we obtain an equation $f(x,x)=a$. We solve it:

Way 1, $f(x,x)=f(x,x+0)=[\mathbf{C}_{2}(f,f_{a})](x,0)=a$, $x=\big\{\mathbf{T}\big[\mathbf{C}_{2}\frac{\mathbf{C}_{2}(f,f_{a})}{0}\big]\big\}^{-1}(a)$, where $\mathbf{C}_{2}\frac{\mathbf{C}_{2}(f,f_{a})}{0}$ is an unary function or we consider it a special binary function.

Way 2, $f(x,x)=f(x,x\times1)=[\mathbf{C}_{2}(f,f_{m})](x,1)=a$, $x=\big\{\mathbf{T}\big[\mathbf{C}_{2}\frac{\mathbf{C}_{2}(f,f_{m})}{1}\big]\big\}^{-1}(a)$, where $\mathbf{C}_{2}\frac{\mathbf{C}_{2}(f,f_{m})}{1}$ is also an unary function.

Anywhere we meet $g[f(x,a),f(x,a)]$ we should converse it to $\{\mathbf{C}_{2}\frac{\mathbf{C}_{2}(g,f_{a})}{0}\}[f$ $(x,a)]=[\mathbf{C}(\mathbf{C}_{2}\frac{\mathbf{C}_{2}(g,f_{a})}{0},f)](x,a)$ or $\{\mathbf{C}_{2}\frac{\mathbf{C}_{2}(g,f_{m})}{1}\}[f(x,a)]=[\mathbf{C}(\mathbf{C}_{2}\frac{\mathbf{C}_{2}(g,f_{m})}{1},f)](x,a)$ then we can reduce the number of f(x,a) from two to one.

\textbf{Example 9} $f_{1}[f_{2}(x,a),f_{3}(x,b)]=c$

Step 1, $f_{1}\rightarrow \mathbf{P}^{3}_{1,3}(f_{1})$,  $f_{2}\rightarrow \mathbf{P}^{3}_{1,2}(f_{2})$, $f_{3}\rightarrow \mathbf{P}^{3}_{1,3}(f_{3})$,

Step 2, $f_{2}(x,a)=\mathbf{P}^{3}_{1,2}(f_{2})(x,a,b)$, really existing variables,x,a.

 \qquad\qquad$f_{3}(x,b)=\mathbf{P}^{3}_{1,3}(f_{3})(x,a,b)$, really existing variables,x,b.

$f_{1}[f_{2}(x,a),f_{3}(x,b)]=\mathbf{P}^{3}_{1,3}[f_{1})[f_{2}(x,a),a,f_{3}(x,b)]$,really existing variables,
$f_{2}(x,a)$,$f_{3}(x,b)$.

$\mathbf{P}^{3}_{1,3}(f_{1})[f_{2}(x,a),a,f_{3}(x,b)]=\mathbf{P}^{3}_{1,3}(f_{1})[P^{3}_{1,2}(f_{2})(x,a,b),a,\mathbf{P}^{3}_{1,3}(f_{3})(x,a,b)]$

$=[\mathbf{C}_{3}\frac{\mathbf{C}_{1}(\mathbf{P}^{3}_{1,3}(f_{1}),\mathbf{P}^{3}_{1,2}(f_{2}))}{\mathbf{P}^{3}_{1,3}(f_{3})}](x,a,b)$

Step3, $x=\mathbf{I}_{1}\{\mathbf{T}[\mathbf{C}_{3}\frac{\mathbf{C}_{1}(\mathbf{P}^{3}_{1,3}(f_{1}),\mathbf{P}^{3}_{1,2}(f_{2}))}{\mathbf{P}^{3}_{1,3}(f_{3})}]\}(c,a,b)$

Special situation,b=a,$f_{1}[f_{2}(x,a),f_{3}(x,a)]=c$,

Use composition of example 3,we have:

$f_{1}[f_{2}(x,a),f_{3}(x,a)]=\mathbf{C}_{3}\Big[C_{3}\frac{\mathbf{C}_{1}(\mathbf{P}^{3}_{1,3}(f_{1}),\mathbf{P}^{3}_{1,2}(f_{2}))}{\mathbf{P}^{3}_{1,3}(f_{3})},
a\Big](x,a)=c$

$\mathbf{C}_{3}\Big[\mathbf{C}_{3}\frac{\mathbf{C}_{1}(\mathbf{P}^{3}_{1,3}(f_{1}),\mathbf{P}^{3}_{1,2}(f_{2}))}{\mathbf{P}^{3}_{1,3}(f_{3})},a\Big]$ is a binary function,a special one of 3 variables.

There are binary functions being unable to be obtained by binary function composition,like $\mathbf{C}_{3}\Big[\mathbf{C}_{3}\frac{\mathbf{C}_{1}(\mathbf{P}^{3}_{1,3}(f_{1}),\mathbf{P}^{3}_{1,2}(f_{2}))}{\mathbf{P}^{3}_{1,3}(f_{3})},a\Big]$.

It needs composition of functions of more variables to obtain them. Generally speaking,there are functions of n variables being unable to be obtained by only composition of functions of n variables. We shown an example of functions of three variables in the example 10.

We known a saying,the shortest distance between two truths in the real domain is through the complex domain.  Now we can say, the shortest distance between the two truths about the functions of less variables is through the functions of more variables.

$x=\mathbf{I}_{1}\Big\{\mathbf{T}\Big[\mathbf{C}_{3}\Big(\mathbf{C}_{3}\frac{\mathbf{C}_{1}(\mathbf{P}^{3}_{1,3}(f_{1}),\mathbf{P}^{3}_{1,2}(f_{2}))}
{\mathbf{P}^{3}_{1,3}(f_{3})},a\Big)\Big]\Big\}(c,a)$

\textbf{Example 10} $f_{1}[f_{2}(x,x_{2},x_{3}),f_{3}(x,x_{4},x_{5}),f_{4}(x,x_{6},x_{7})]=a$

This is an equation constructed by functions of three variables.

Step 1: $f_{1}\rightarrow \mathbf{P}^{7}_{1,4,6}(f_{1})$, $f_{2}\rightarrow \mathbf{P}^{7}_{1,2,3}(f_{2})$, $f_{3}\rightarrow \mathbf{P}^{7}_{1,4,5}(f_{3})$, $f_{4}\rightarrow \mathbf{P}^{7}_{1,6,7}(f_{4})$

Step 2:   $f_{2}(x,x_{2},x_{3})=\mathbf{P}^{7}_{1,2,3}(f_{2})(x,x_{2},x_{3},x_{4},x_{5},x_{6},x_{7})$,really existing variables,x, $x_{2}$ and $x_{3}$.

 $f_{3}(x,x_{4},x_{5})=\mathbf{P}^{7}_{1,4,5}(f_{2})(x,x_{2},x_{3},x_{4},x_{5},x_{6},x_{7})$,really existing variables,x,$x_{4}$,$x_{5}$.

 $f_{4}(x,x_{6},x_{7})=\mathbf{P}^{7}_{1,6,7}(f_{2})(x,x_{2},x_{3},x_{4},x_{5},x_{6},x_{7})$,really existing variables,x,$x_{6}$,$x_{7}$.

$f_{1}[f_{2}(x,x_{2},x_{3}),f_{3}(x,x_{4},x_{5}),f_{4}(x,x_{6},x_{7})]$

$=\mathbf{P}^{7}_{1,4,6}(f_{1})\Big\{\underbrace{f_{2}(x,x_{2},x_{3})}_{1th.variable},x_{2},x_{3},\underbrace{f_{3}(x,x_{4},x_{5})}_{4th.variable},
x_{5},x_{6},\underbrace{f_{4}(x,x_{6},x_{7})}_{6th.variable}\Big\}$

$=\mathbf{P}^{7}_{1,4,6}(f_{1})\Big\{\underbrace{\mathbf{P}^{7}_{1,2,3}(f_{2})(x,x_{2},x_{3},x_{4},x_{5},x_{6},x_{7})}_{1th.variable},x_{2},x_{3},\underbrace{
\mathbf{P}^{7}_{1,4,5}(f_{2})(x,x_{2},x_{3},x_{4},x_{5},x_{6},x_{7})}_{4th.variable},$

$x_{5},x_{6},\underbrace{\mathbf{P}^{7}_{1,6,7}(f_{2})(x,x_{2},x_{3},x_{4},x_{5},x_{6},x_{7})}_{6th.variable}\Big\}$

$=\mathbf{C}_{6}\Big[\mathbf{C}_{4}\frac{\mathbf{C}_{1}[\mathbf{P}^{7}_{1,4,6}(f_{1}),\mathbf{P}^{7}_{1,2,3}(f_{2})]}{\mathbf{P}^{7}_{1,4,5}(f_{3})},
\mathbf{P}^{7}_{1,6,7}(f_{4})\Big](x,x_{2},x_{3},x_{4},x_{5},x_{6},x_{7})$

Step 3:

$x=\mathbf{I}_{1}\Big\{\mathbf{T}\Big[\mathbf{C}_{6}\Big(\mathbf{C}_{4}\frac{\mathbf{C}_{1}[\mathbf{P}^{7}_{1,4,6}(f_{1}),\mathbf{P}^{7}_{1,2,3}
(f_{2})]}{\mathbf{P}^{7}_{1,4,5}(f_{3})},\mathbf{P}^{7}_{1,6,7}(f_{4})\Big)\Big]\Big\}(a,x_{2},x_{3},x_{4},x_{5},x_{6},x_{7})$

Special situation,$x_{4}=x_{6}=x_{2}$,$x_{5}=x_{7}=x_{3}$,

$f_{1}[f_{2}(x,x_{2},x_{3}),f_{3}(x,x_{2},x_{3}),f_{4}(x,x_{2},x_{3})]=a$

Use composition of example 3 four times,we have:

$f_{1}[f_{2}(x,x_{2},x_{3}),f_{3}(x,x_{2},x_{3}),f_{4}(x,x_{2},x_{3})]$

$=\mathbf{C}_{6}\big[\mathbf{C}_{4}\frac{\mathbf{C}_{1}[\mathbf{P}^{7}_{1,4,6}(f_{1}),\mathbf{P}^{7}_{1,2,3}(f_{2})]}{\mathbf{P}^{7}_{1,4,5}(f_{3})},
\mathbf{P}^{7}_{1,6,7}(f_{4})\big](x,x_{2},x_{3},x_{2},x_{3},x_{2},x_{3})$

$=\big\{\mathbf{C}_{7}\frac{\mathbf{C}_{6}\big[\mathbf{C}_{4}\frac{\mathbf{C}_{1}[\mathbf{P}^{7}_{1,4,6}(f_{1}),
\mathbf{P}^{7}_{1,2,3}(f_{2})]}{\mathbf{P}^{7}_{1,4,5}(f_{3})},\mathbf{P}^{7}_{1,6,7}(f_{4})\big]}{x_{3}}\big\}
(x,x_{2},x_{3},x_{2},x_{3},x_{2})$

$=\bigg(\mathbf{C}_{6}\big\{\mathbf{C}_{7}\frac{\mathbf{C}_{6}\big[\mathbf{C}_{4}\frac{\mathbf{C}_{1}[\mathbf{P}^{7}_{1,4,6}(f_{1}),\mathbf{P}^{7}_{1,2,3}
(f_{2})]}{\mathbf{P}^{7}_{1,4,5}(f_{3})},\mathbf{P}^{7}_{1,6,7}(f_{4})\big]}{x_{3}}\big\},x_{2}\bigg)(x,x_{2},x_{3},x_{2},x_{3})$

$=\bigg[\frac{\mathbf{C}_{5}\bigg(\mathbf{C}_{6}\big\{\mathbf{C}_{7}\frac{\mathbf{C}_{6}\big[\mathbf{C}_{4}\frac{\mathbf{C}_{1}[\mathbf{P}^{7}_{1,4,6}(f_{1}),
\mathbf{P}^{7}_{1,2,3}(f_{2})]}{\mathbf{P}^{7}_{1,4,5}(f_{3})},\mathbf{P}^{7}_{1,6,7}(f_{4})\big]}{x_{3}}\big\},x_{2}\bigg)}{x_{3}}\bigg](x,x_{2},x_{3},x_{2})$

$=\bigg\{\mathbf{C}_{4}\bigg[\frac{\mathbf{C}_{5}\bigg(\mathbf{C}_{6}\big\{\mathbf{C}_{7}\frac{\mathbf{C}_{6}\big[\mathbf{C}_{4}\frac{\mathbf{C}_{1}
[\mathbf{P}^{7}_{1,4,6}(f_{1}),\mathbf{P}^{7}_{1,2,3}(f_{2})]}{\mathbf{P}^{7}_{1,4,5}(f_{3})},\mathbf{P}^{7}_{1,6,7}(f_{4})\big]}{x_{3}}\big\},x_{2}
\bigg)}{x_{3}}\bigg],x_{2}\bigg\}(x,x_{2},x_{3})$

$\bigg\{\mathbf{C}_{4}\bigg[\frac{\mathbf{C}_{5}\bigg(\mathbf{C}_{6}\big\{\mathbf{C}_{7}\frac{\mathbf{C}_{6}\big[\mathbf{C}_{4}\frac{\mathbf{C}_{1}
[\mathbf{P}^{7}_{1,4,6}(f_{1}),\mathbf{P}^{7}_{1,2,3}(f_{2})]}{\mathbf{P}^{7}_{1,4,5}(f_{3})},\mathbf{P}^{7}_{1,6,7}(f_{4})\big]}{x_{3}}\big\},x_{2}
\bigg)}{x_{3}}\bigg],x_{2}\bigg\}$is a function of 3 variables,a special one of 7 variables.

$x=\mathbf{I}_{1}\Bigg(\mathbf{T}\bigg\{\mathbf{C}_{4}\bigg[\frac{\mathbf{C}_{5}\bigg(\mathbf{C}_{6}\big\{\mathbf{C}_{7}\frac{\mathbf{C}_{6}\big[\mathbf{C}_{4}\frac{\mathbf{C}_{1}
[\mathbf{P}^{7}_{1,4,6}(f_{1}),\mathbf{P}^{7}_{1,2,3}(f_{2})]}{\mathbf{P}^{7}_{1,4,5}(f_{3})},\mathbf{P}^{7}_{1,6,7}(f_{4})\big]}{x_{3}}\big\},x_{2}
\bigg)}{x_{3}}\bigg],x_{2}\bigg\}(x,x_{2},x_{3})\Bigg)(a,x_{2},x_{3})$

\textbf{Example 11} $f_{1}[f_{2}(x,a),b]+f_{3}[f_{4}(x,c),d]=e$

Step 1, $f_{1}\rightarrow \mathbf{P}^{5}_{1,3}(f_{1})$,  $f_{2}\rightarrow \mathbf{P}^{5}_{1,2}(f_{2})$, $f_{3}\rightarrow \mathbf{P}^{5}_{4,5}(f_{3})$, $f_{4}\rightarrow \mathbf{P}^{5}_{1,4}(f_{4})$, $f_{a}\rightarrow \mathbf{P}^{5}_{1,4}(f_{a})$

Step 2, $ f_{2}(x,a)=\mathbf{P}^{5}_{1,2}(f_{2})(x,a,b,c,d)$,really existing variables,x,a.

$f_{1}[f_{2}(x,a),b]=\mathbf{P}^{5}_{1,3}(f_{1})[f_{2}(x,a),a,b,c,d]$,really existing variables,$f_{2}(x,a)$,b.

$f_{1}[f_{2}(x,a),b]=\mathbf{P}^{5}_{1,3}(f_{1})[\mathbf{P}^{5}_{1,2}(f_{2})(x,a,b,c,d),a,b,c,d]=\mathbf{C}_{1}[\mathbf{P}^{5}_{1,3}(f_{1}),\mathbf{P}^{5}_{1,2}(f_{2})](x,a,b,c,d)$,

 $ f_{4}(x,c)=\mathbf{P}^{5}_{1,4}(f_{4})(x,a,b,c,d)$, really existing variables, x,c.

$f_{3}[f_{4}(x,c),d]=\mathbf{P}^{5}_{4,5}(f_{3})[f_{4}(x,c),a,b,c,d]$, really existinging variables,$f_{4}(x,c)$,c.

$f_{3}[f_{4}(x,c),d]=\mathbf{P}^{5}_{4,5}(f_{3})[\mathbf{P}^{5}_{1,4}(f_{4})(x,a,b,c,d),a,b,c,d]=\mathbf{C}_{1}[\mathbf{P}^{5}_{4,5}(f_{3}),\mathbf{P}^{5}_{1,4}(f_{4})](x,a,b,c,d)$,

 $f_{1}[f_{2}(x,a),b]+f_{3}[f_{4}(x,c),d]=f_{a}\{f_{1}[f_{2}(x,a),b],f_{3}[f_{4}(x,c),d]\}$

 $=\mathbf{P}^{5}_{1,4}(f_{a})\{\underbrace{f_{1}[f_{2}(x,a),b]}_{1th.variable},a,b,\underbrace{f_{3}[f_{4}(x,c),d]}_{4th.variable},d\}$,

 really existing variables, $f_{1}[f_{2}(x,a),b]$,$f_{3}[f_{4}(x,c),d]$.

 $\mathbf{P}^{5}_{1,4}(f_{a})\{f_{1}[f_{2}(x,a),b],a,b,f_{3}[f_{4}(x,c),d],d\}$

 $=\mathbf{P}^{5}_{1,4}(f_{a})\{\underbrace{\mathbf{C}_{1}[\mathbf{P}^{5}_{1,3}(f_{1}),\mathbf{P}^{5}_{1,2}(f_{2})](x,a,b,c,d)}_{1th.variable},a,b,
 \underbrace{\mathbf{C}_{1}[\mathbf{P}^{5}_{4,5}(f_{3}),\mathbf{P}^{5}_{1,4}(f_{4})](x,a,b,c,d)}_{4th.variable},d\}$

$=\mathbf{C}_{4}\Big[\Big(\mathbf{C}_{1}\frac{\mathbf{P}^{5}_{1,4}(f_{a})}{\{\mathbf{C}_{1}[\mathbf{P}^{5}_{1,3}(f_{1}),\mathbf{P}^{5}_{1,2}
(f_{2})]\}}\Big),\mathbf{C}_{4}[\mathbf{P}^{5}_{4,5}(f_{3}),\mathbf{P}^{5}_{1,4}(f_{4})]\Big](x,a,b,c,d)$

Step 3,

$x=\Bigg[\mathbf{I}_{1}\Bigg(\mathbf{T}\Big\{\mathbf{C}_{4}\Big[\Big(\mathbf{C}_{1}\frac{\mathbf{P}^{5}_{1,4}(f_{a})}{\{\mathbf{C}_{1}[\mathbf{P}^{5}_{1,3}
(f_{1}),\mathbf{P}^{5}_{1,2}(f_{2})]\}}\Big),\mathbf{C}_{4}[\mathbf{P}^{5}_{4,5}(f_{3}),\mathbf{P}^{5}_{1,4}(f_{4})]\Big]\Big\}\Bigg)\Bigg](e,a,b,c,d)$

\textbf{Example 12} $ax^{2}+bx=c$

We omit the procedure and give the result directly for this example and the next one.

$x=\Bigg(\mathbf{I}_{2}\Big\{\mathbf{T}\Big[\mathbf{C}_{2}\Big(\mathbf{C}_{1}\frac{\mathbf{P}^{4}_{1,4}(f_{a})}{\mathbf{C}_{2}[\mathbf{P}^{4}_{1,2}(f_{m}),
\mathbf{P}^{4}_{2,3}(f_{p})]},\mathbf{P}^{4}_{2,4}(f_{m})\Big)\Big]\Big\}\Bigg)(a,c,2,b)$

Surely there is binary solution for it but we do not discuss here.

\textbf{Example 13} $x^{7}+ax^{3}+bx^{2}+cx=-1$

This equation was mentioned in famous Hilbert's 13th problem.

$x=\mathbf{I}_{1}\bigg[\mathbf{T}\bigg(\mathbf{C}_{7}\frac{\bigg[\mathbf{C}_{1}\bigg(\mathbf{P}^{7}_{1,7}(f_{a}),\mathbf{C}_{5}\frac{\Big[C_{1}\Big(
\mathbf{P}^{7}_{1,5}(f_{a}),\mathbf{C}_{3}\frac{\mathbf{C}_{1}[P^{7}_{1,3}(f_{a}),\mathbf{P}^{7}_{1,2}(f_{p})]}{\mathbf{C}_{4}[P^{7}_{3,4}(f_{m}),
\mathbf{P}^{7}_{1,4}(f_{p})]}\Big)\Big]}{\mathbf{C}_{6}[\mathbf{P}^{7}_{5,6}(f_{m}),\mathbf{P}^{7}_{1,6}(f_{p})]}\bigg)\bigg]}{\mathbf{P}^{7}_{7,1}
(f_{m})}\bigg)\bigg](-1,7,a,3,b,2,c)$

\textbf{Example 14} $f_{1}(x,a)+f_{2}(x,a)=b$

 $f_{1}(x,a)+f_{2}(x,a)=b\Longrightarrow (f_{1}+f_{2})(x,a)=b$

$f_{1}+f_{2}$ is called the superposition of $f_{1}$ and $f_{2}$.

$x=\{\mathbf{I}_{1}[\mathbf{T}(f_{1}+f_{2})]\}(b,a),j=1$. This is the simplest form of the solution.

Can we use the inverse symbolic $\mathbf{I}_{i}$ in the expression of the solution? We say we must use it! For the equation $R_{1}$(x,a,b,c)=d,x=$R_{2}$ (d,a,b,c). $R_{2}$ is a new function or a new relation. How do we express it? By an algebraic function of a,b c and d? It does not exist for a general equation.
We want to give explicit solution for a general equation we have to create new form of functions or relations.  However, it is impossible to name every new function or new relation. Pay attention to obtaining these new functions or new relations by ones we have had. It is a good idea to express them by the old function or old relation, $R_{2}=\mathbf{I}_{i}(R_{1})$. Why could not we do like this? There are seven operations for numbers,addition,substraction, multiplication,division, power, root and logarithm. However there are only four operations for functions, $\mathbf{C}_{i}$, $\mathbf{P}^{n}_{j_{1},j_{2},\cdots,j_{m}}$, $\mathbf{I}_{i}$ and $\mathbf{T}$. Root and logarithm may be reject in the history. However they were accepted finally and actually they had been used up to now. Likewise, $\mathbf{C}_{i}$, $\mathbf{P}^{n}_{j_{1},j_{2},\cdots,j_{m}}$, $\mathbf{I}_{i}$ and $\mathbf{T}$ will be accepted too.

Never mind how complex form the solutions have. We can omitted symbolics $\mathbf{I}_{i}$ and $\mathbf{T}$ existing once and being always in the left of the solution. Symbolic $\mathbf{P}^{n}_{j1,j2,\cdots,jm}$ is always in the front of a binary function. Symbolic $C_{i}$ are always be operated from the right to the left.
 $\sqrt{2}$ and $\sqrt{5}$  are numbers,$\sqrt{2}+\sqrt{5}$ is also a number with a complex form. Simple functions $f_{1}$ and $f_{2}$ are functions, $C_{2}(f_{1},f_{2})$,the function obtained by a simple composition and  $\bigg\{\mathbf{C}_{4}\bigg[\frac{\mathbf{C}_{5}\bigg(\mathbf{C}_{6}\big\{\mathbf{C}_{7}\frac{\mathbf{C}_{6}\big[\mathbf{C}_{4}\frac{\mathbf{C}_{1}
[\mathbf{P}^{7}_{1,4,6}(f_{1}),\mathbf{P}^{7}_{1,2,3}(f_{2})]}{\mathbf{P}^{7}_{1,4,5}(f_{3})},\mathbf{P}^{7}_{1,6,7}(f_{4})\big]}{x_{3}}\big\},x_{2}
\bigg)}{x_{3}}\bigg],x_{2}\bigg\}(x,x_{2},x_{3})$, the function created by several compositions are also functions.

We promised to give ways to express the multivariate solution for a general equation. Here we did and shown several type of multivariate solutions. Are these expressions different from useless symbolics f and g? These expressions remains all elements of the equations!

{\bf 3.2 Types of multivariate functions or of multivariate relations}

There are four kinds of multivariate functions or of multivariate relations.

Type 1, it and its inverse can be expressed by binary function composition without symbolic of inverse.

Type 2, it  can be expressed by binary function composition without symbolic of inverse otherwise its inverse can not be expressed by binary function composition without symbolic of inverse,such as $(xf_{1}a)f_{3}(xf_{2}b)$.

Type 3, it  can not be expressed by binary function composition without symbolic of inverse otherwise its inverse can be expressed by binary function composition without symbolic of inverse.

Type 4, it and its inverse can not be expressed by binary function composition without symbolic of inverse. For example, $f _{i},i=1\sim5$ are general binary functions, $\mathbf{I}_{1}\Big[\mathbf{C}_{3}\Big(\mathbf{C}_{3}\frac{\mathbf{C}_{1}(\mathbf{P}^{3}_{1,3}(f_{1}),\mathbf{P}^{3}_{1,2}(f_{2}))}
{\mathbf{P}^{3}_{1,3}(f_{3})},a\Big)\Big]+\mathbf{C}_{3}$ $\Big(\mathbf{C}_{3}\frac{\mathbf{C}_{1}(\mathbf{P}^{3}_{1,3}(f_{1}),\mathbf{P}^{3}_{1,2}(f_{2}))}
{\mathbf{P}^{3}_{1,3}(f_{3})},a\Big)$ must be in type 4.

The inverse of a function or relation of type 1 must be in type 1.The inverse of a function or relation of type 2 must be in Type 3.The inverse of a function or relation of type 3 must be in type 2.The inverse of a function or relation of type 4 must be in type 4. We can express it like this: $1\Longrightarrow1$,\qquad$2\Longrightarrow3$,\qquad$3\Longrightarrow2$,\qquad$4\Longrightarrow4$

Which are the main roles when we study equations? Binary ones or multivariate ones? Binary ones are problems. So many equations are created
by them. However we can not use them to express the solutions of these equations. Multivariate ones are answers. We can use them to  express the solutions
of any equations. We ask again,which are the main roles? Binary ones or multivariate ones? Please forget binary ones. They are only transitional roles. Multivariate ones are the main roles. Kingdom of equations will be perfect if and only if it includes all multivariate ones.

\section{Binary solutions for special equations}

There are many equations with binary solution, such as example 8, example 9 and example 14, polynomial algebraic equations(n$\leq$6), etc. Our target is to find more equations with binary solution. In this section we will show an important kind of equations with binary solution.

$\mathbf{Definition~4.1}$ Two binary functions $f_{1}$ and $f_{2}$ are $\mathbf{consistent}$ if they meet

$f_{1}[f_{2}(x_{1},x_{2}),x_{3}]=f_{2}[f_{1}(x_{1},x_{3}),x_{2}]$.

 $\mathbf{Definition~4.2}$ A set of functions F is called  $\mathbf{cognate}$ if $\forall f_{1},f_{2}\in F$ are $\mathbf{consistent}$.

\textbf{Example 15}  $F:=\{\beta|\beta=C_{2}[f_{p},c]\}$,F is $\mathbf{consistent}$. See \textbf{Example 2}, Where $\beta$ is obtained from cutting the space surface of $f_{p}$ with plane parallel to the XoZ coordinate surface.

If $f_{1},f_{3}\in F$, $f_{1}$ and $f_{3}$ are consistent then $(xf_{1}a)f_{3}(xf_{2}b)=c$ has binary solution.

$(xf_{1}a)f_{3}(xf_{2}b)=[xf_{3}(xf_{2}b)]f_{1}a)=\{[\mathbf{C}_{2}(f_{3},f_{2})](x,b)\}f_{1}a=c$

$[\mathbf{C}_{2}(f_{3},f_{2})](x,b)=[\mathbf{I}_{1}(f_{1})](c,a)$

$x=\{\mathbf{I}_{1}[\mathbf{C}_{2}(f_{3},f_{2})]\}\bigg(\{[\mathbf{I}_{1}(f_{1})](c,a)\},b\bigg)$

If $f_{2}$ and $f_{3}$ are consistent then $(xf_{1}a)f_{3}(xf_{2}b)=c$ has binary solution too.

\section{Express multivariate solutions by superposition of unary functions}

We mentioned equations with binary solutions above section. These  binary solutions are perfect for they being constructed by binary function with good properties. May be these binary function are analytic functions or smooth functions. If we loose the condition to continuous functions then there are binary solutions or unary function for all of equations by following theorem.

{\bf 5.1 Expressing a continuous multivariate function by superposition of unary functions}

 {\bf Theorem~5.1~~Kolmogorov's Superposition Theorem} Let $f:[0,1]^{n}$ $\longrightarrow R$ be an arbitrary multivariate continuous function.Then it has the representation.

\begin{equation}\label{eq:eps}
f(x_{1},x_{2},\cdots,x_{n})=\sum _{q=0}^{2n}f_{q}\big[\sum _{p=1}^{n}g_{q,p}(x_{p})\big]
\end{equation}

 Here a set  of inner functions $g_{q,p}$ is not unique but any selected set will be independent to $f$ then they are called remaining functions.

 {\bf Definition~5.1 Inner set:} A set of inner functions $g_{q,p}$ is called a inner set.

{\bf Definition~5.2 Multivariate function decomposition } To express the relation between $f_{q}$ and f, we introduce decomposition operator $D_{q}$��

\begin{equation}\label{eq:eps}
   f_{q}=\mathbf{D}_{q}(f),(0\leq q\leq 2n)
\end{equation}

$f$ can be expressed too:

\begin{equation}\label{eq:eps}
f(x_{1},x_{2},\cdots,x_{n})=\sum _{q=0}^{2n}\mathbf{D}_{q}(f)\big[\sum _{p=1}^{n}g_{q,p}(x_{p})\big]
\end{equation}

 This expression is called multivariate function decomposition or superposition of unary functions of it.

The binary solution for example 13 should be:

$x=\bigg\{\mathbf{I}_{1}\bigg[\mathbf{T}\bigg(\mathbf{C}_{7}\frac{\bigg[\mathbf{C}_{1}\bigg(\mathbf{P}^{7}_{1,7}(f_{a}),\mathbf{C}_{5}\frac{\Big[C_{1}\Big(
\mathbf{P}^{7}_{1,5}(f_{a}),\mathbf{C}_{3}\frac{\mathbf{C}_{1}[P^{7}_{1,3}(f_{a}),\mathbf{P}^{7}_{1,2}(f_{p})]}{\mathbf{C}_{4}[P^{7}_{3,4}(f_{m}),
\mathbf{P}^{7}_{1,4}(f_{p})]}\Big)\Big]}{\mathbf{C}_{6}[\mathbf{P}^{7}_{5,6}(f_{m}),\mathbf{P}^{7}_{1,6}(f_{p})]}\bigg)\bigg]}{\mathbf{P}^{7}_{7,1}
(f_{m})}\bigg)\bigg]\bigg\}(-1,7,a,3,b,2,c)$

$$=\sum _{q=0}^{14}\mathbf{D}_{q}\bigg\{\mathbf{I}_{1}\bigg[\mathbf{T}\bigg(\mathbf{C}_{7}\frac{\bigg[\mathbf{C}_{1}\bigg(\mathbf{P}^{7}_{1,7}(f_{a}),\mathbf{C}_{5}\frac{\Big[C_{1}\Big(
\mathbf{P}^{7}_{1,5}(f_{a}),\mathbf{C}_{3}\frac{\mathbf{C}_{1}[P^{7}_{1,3}(f_{a}),\mathbf{P}^{7}_{1,2}(f_{p})]}{\mathbf{C}_{4}[P^{7}_{3,4}(f_{m}),
\mathbf{P}^{7}_{1,4}(f_{p})]}\Big)\Big]}{\mathbf{C}_{6}[\mathbf{P}^{7}_{5,6}(f_{m}),\mathbf{P}^{7}_{1,6}(f_{p})]}\bigg)\bigg]}{\mathbf{P}^{7}_{7,1}
(f_{m})}\bigg)\bigg]\bigg\}$$

$$\bigg[g_{q,1}(-1)+g_{q,2}(7)+g_{q,3}(a)+g_{q,4}(3)+g_{q,5}(b)+g_{q,6}(2)+g_{q,7}(c)\bigg],j=1\sim 7$$

$$Note,\mathbf{D}_{q}\bigg\{\mathbf{I}_{1}\bigg[\mathbf{T}\bigg(\mathbf{C}_{7}\frac{\bigg[\mathbf{C}_{1}\bigg(\mathbf{P}^{7}_{1,7}(f_{a}),\mathbf{C}_{5}\frac{\Big[C_{1}\Big(
\mathbf{P}^{7}_{1,5}(f_{a}),\mathbf{C}_{3}\frac{\mathbf{C}_{1}[P^{7}_{1,3}(f_{a}),\mathbf{P}^{7}_{1,2}(f_{p})]}{\mathbf{C}_{4}[P^{7}_{3,4}(f_{m}),
\mathbf{P}^{7}_{1,4}(f_{p})]}\Big)\Big]}{\mathbf{C}_{6}[\mathbf{P}^{7}_{5,6}(f_{m}),\mathbf{P}^{7}_{1,6}(f_{p})]}\bigg)\bigg]}{\mathbf{P}^{7}_{7,1}
(f_{m})}\bigg)\bigg]\bigg\}$$ are just only unary functions.

This expression of the solution contains only known parameters,known functions or known operations, remaining functions and function promotion,
composition,inverse,decomposition.There are infinite inner sets. However there must a special inner set,like (0,0,1),(0,1,0),(1,0,0) for 3 dimension space, it has good characteristics that other inner sets must lack. What good characteristics it has and how to find this special inner set are also open problems for us.

   Must we  define decomposition? We must inject new elements if they are necessary. We can not do anything if we reject these new elements.We say that we could not
study trigonometrical function if we had not accepted sin,cos,tan,cotan. Abel and Jacobi  introduced the elliptical function to express the inverse of elliptical integral then they made a big progress in their study. Meanwhile Legendre was confused about this key point then he wasted several years. In the same reason we will miss the chance of expressing the solution of a general equation by unary functions if we refuse to accept decomposition.

     Unary functions in Kolmogorov's superposition theorem are continuous so we can obtain only solutions expressed by unary continuous function composition.
These unary continuous functions lack good property for they being very steep in some ranges. We lose something when we obtain solutions expressed by binary functions or unary functions. This is reasonable. There is no any necessary for  multivariate functions to existing if arbitrary  multivariate functions can be expressed by good binary functions, such as binary analytic functions or binary algebraic functions. Is God  generous or mean? He gave equations binary solutions, unluckily, not perfect ones.

\section{Discrete equations}

    Discrete equation are a special kind of equations. We have the theorem that an arbitrary multivariate discrete function can be expressed by binary discrete
functions or by unary discrete functions. Many problems became simple if they are about discrete equations. You can see this point from our results mentioned in an other paper 'DISCRETE EQUATIONS'.

\section{Discussion and expectation}

Solving an equation is reducing the number of unknown X to one then placing this alone X to the left of '=' and placing all the known items to the right of '='. The multivariate solution given in this paper meets this requirement. Achieving this is a huge breakthrough.
Multivariate function composition and inverse multivariate function were introduced specially by author to solve transcendental equations. Human can reject these new concepts.  However they will be helpless when they face transcendental equations. They had to admit to having been completely defeated by the equations.
Must we solve transcendental equations by exact formula? Structural mathematics, topology and algebra have been the main role of
mathematics for more than one hundred years since Hilbert's solving Gordan problem. Meanwhile the method of classical mathematics is nearly forgotten
by mathematicians. I think this is a serious imbalance. However resources in topology and algebra will be dried up in the near future and mathematics must go back to the path of classical mathematics for there being so mineral in it.

 Mathematics has two legs ,the constructive method and the existence method. Both legs should be used interchangeably.  Before Hilbert's solving Gordan problem the constructive method played a major role, mathematics took the left leg.  Since then,  existence method emerged and became the ruler of mathematics,mathematics has moved its right leg. It should be the left leg by now.

We shown perfect way for multivariate solutions of algebraic equations and of operator equations. However there are difficult tasks to find more good binary solutions for special equations and find good way to express multivariate solutions by unary functions. Studying equations is a huge and complicated subject. This paper gave only the basic concepts and basic methods of our theory and described only the outline of it.

We show readers  a ladder road along the cliff to the ground high in the clouds and want to make readers optimistic. But there is a long distance to go if we want to gallop through the vast territory of equations kingdom. Author just make a start about this topic and sincerely hope other ones follow up. No matter how difficult this topic is and whether or not mathematicians snub it, logic will become a decisive force, and the kingdom of equations will surely have a new world.



\qquad \qquad\qquad \qquad\qquad \qquad References 

\bibliographystyle{elsarticle-harv}

\bibliography{<your-bib-database>}

[1] A. G.Vitushkin, On Hilbert's thirteenth problem and related
questions, Russian Math. Surveys 59:1 11-25 (2004)

[2] A.N.Kolmogorov, On the representation of continuous functions of
several variables by superpositions of continuous functions of one
variable and addition, Dokl.Akad.Nauk SSSR 114 (1957), 953-956;
English transl., Amer.Math. Soc.Transl. (2) 28 (1963), 55-59.

[3] D.Hilbert, Mathematical Problems, Bull.Amer.Math.
Soc.8(1902),461-462.







\end{document}